\documentclass[a4paper]{jpconf}

\usepackage{graphicx}
\usepackage{amsfonts}
\usepackage{amsbsy}
\usepackage{pstricks}

\def\ds{\displaystyle}
\def\ov{\overline}
\def\x{x}
\def\y{y}
\def\antip{\textsf{S}}
\def\la{\langle}
\def\ra{\rangle}
\def\lb{[\![}
\def\rb{]\!]}
\def\ch{{\rm ch}\,}
\def\Z{{\mathbb Z}}
\def\C{{\mathbb C}}
\def\cal{\mathcal}
\def\qed{\hfill\fullsquare}

\newtheorem{thm}{Theorem}[section]

\begin{document}
\title{Hopf algebras and characters of classical groups}

\author{Ronald C King$^1$, Bertfried Fauser$^2$, Peter D Jarvis$^3$}

\address{$^1$ School of Mathematics, University of Southampton, Southampton, SO17 1BJ, England}
\address{$^2$ Max Planck Institut F\"{u}r Mathematik, Inselstrasse 22-26, D-04103 Leipzig, Germany}
\address{$^3$ School of Mathematics and Physics, University of Tasmania, GPO Box 252-21, 7001 Hobart, Tasmania, Australia}
\ead{R.C.King@soton.ac.uk, fauser@mis.mpg.de, Peter.Jarvis@utas.edu.au}

\begin{abstract}
Schur functions provide an integral basis of the ring of symmetric functions.
It is shown that this ring has a natural Hopf algebra structure by  
identifying the appropriate product, coproduct, unit, counit and antipode, 
and their properties. 
Characters of covariant tensor irreducible representations of the
classical groups $GL(n)$, $O(n)$ and $Sp(n)$ are then expressed
in terms of Schur functions, and the Hopf algebra is exploited in the 
determination of group-subgroup branching rules and the decomposition
of tensor products. The analysis is carried out in terms of
$n$-independent universal characters. The corresponding rings,
$CharGL$, $CharO$ and $CharSp$, of universal characters
each have their own natural Hopf algebra structure. 
The appropriate product, coproduct, unit, counit and antipode are
identified in each case.
\end{abstract}


\section{Introduction}
The aim here is to provide a uniform setting for dealing with
characters of finite-dimensional irreducible representations of the 
classical groups $GL(n)$, $O(n)$ and $Sp(n)$. More specifically,
the intention is to determine certain group-subgroup branching rules
and formulae for the decomposition of tensor products of irreducible 
representations of all of these groups by using Schur functions and  
the Hopf algebra of these symmetric functions, as described most 
recently in~\cite{FJKW,FJK}

\section{The Hopf algebra of symmetric functions}

Let $x=(x_1,x_2,\ldots,x_n)$ be a sequence of $n$ indeterminates
and let $\Lambda^{(n)}=\Z[\x]^{S_n}$ be the ring of polynomial symmetric 
functions of the indeterminates $x_1,x_2,\ldots,x_n$. This
ring may be graded by the total degree, $d$, of such polynomials,
so that we may write $\Lambda^{(n)}=\oplus_d\,\Lambda^{(n)}_d$. 
An integral basis of $\Lambda^{(n)}_d$ is provided by the
Schur functions~\cite{Li50,Mac} 
\begin{equation}
  s_\lambda(x)=\frac{\left|\,x_i^{\lambda_j+n-j}\,\right|}{\left|\,x_i^{n-j}\,\right|}\,,
\label{eq-sfn}
\end{equation}
specified by partitions 
$\lambda=(\lambda_1,\lambda_2,\ldots,\lambda_n)$ 
of weight $|\lambda|=\lambda_1+\lambda_2+\cdots+\lambda_n=d$
and length $\ell(\lambda)=p\leq n$, so that 
$\lambda_1\geq\lambda_2\geq\cdots\geq\lambda_p>0$
with $\lambda_i=0$ for any $i>p$. 

Within the ring $\Lambda^{(n)}$, products of 
Schur functions decompose as follows:
\begin{equation}
s_\lambda(\x)\ s_\mu(\x)
    = \sum_\nu\ c_{\lambda\mu}^\nu\ s_\nu(\x)\,,
\label{eq-LRcoeff}    
\end{equation}
where the coefficients $c_{\lambda\mu}^\nu$ are known as Littlewood-Richardson
coefficients. They are all non-negative integers and may be evaluated by means of the 
Littlewood-Richardson rule~\cite{Li50,Mac}.

If we move to a sequence $x=(x_1,x_2,\ldots)$ of countably many independent
indeterminates, then for all partitions $\lambda$ of weight $|\lambda|=d$ with $d$ finite,
there exists a universal Schur function $s_\lambda(x)$ of $x=(x_1,x_2,\ldots)$~\cite{Mac}
such that for all finite $n$ we have 
$s_\lambda(x_1,x_2,\ldots,x_n,0,0,\ldots)=s_\lambda(x_1,x_2,\ldots,x_n)\in \Lambda^{(n)}_d$.
This stability property enables us to define the ring $\Lambda$ of symmetric
functions. This is the ring generated by $s_\lambda(x)$ for all partitions $\lambda$.
Within this ring $\Lambda$, the multiplication rule is again given by (\ref{eq-LRcoeff}),
and still within $\Lambda$, skew Schur functions~\cite{Li50,Mac} are defined by:
\begin{equation}
s_{\nu/\lambda}(x)=\sum_{\mu}\, c_{\lambda\mu}^\nu\, s_\mu(x)\,.
\label{eq-skew}
\end{equation}

Each partition $\lambda$ defines a Young diagram $F^\lambda$ whose 
successive rows lengths are the parts of $\lambda$, and whose 
successive column lengths are the parts of the conjugate 
partition, denoted here by $\lambda'$. 
Then for all $\x=(x_1,x_2,\ldots)$ and $\y=(y_1,y_2,\ldots)$
Cauchy's formula and its inverse take the form~\cite{Li50,Mac}:
\begin{eqnarray}
&&J(x,y)=\prod_{i,a}\ (1-x_iy_a)^{-1}=\sum_\lambda\,s_\lambda(\x)\,s_\lambda(\y);
\label{eq-J}\\
&&I(x,y)=\prod_{i,a}\ (1-x_iy_a)=\sum_\lambda\,(-1)^{|\lambda|}\,s_\lambda(\x)\,s_{\lambda'}(\y).
\label{eq-I}
\end{eqnarray}
It follows that
\begin{eqnarray*}
J(x,y)\ I(x,y)
&=&\sum_\lambda\,s_\lambda(\x)\,s_\lambda(\y)\ \sum_\mu\,(-1)^{|\mu|}\,s_\mu(\x)\,s_{\mu'}(\y)\\
&=&\sum_{\lambda,\mu,\nu} c_{\lambda\mu}^\nu\,s_\nu(\x)\ (-1)^{|\mu|}
\,s_\lambda(\y)\,s_{\mu'}(\y)
=\sum_{\mu,\nu} s_\nu(\x)\ (-1)^{|\mu|}\,s_{\nu/\mu}(\y)\,s_{\mu'}(\y)\,.
\end{eqnarray*}
However $J(x,y)\ I(x,y)=1=s_0(\x)\,s_0(\y)$, so that by
comparing coefficients of $s_\nu(x)$ we have the Schur function identity:
\begin{equation}
\sum_\mu\, (-1)^{|\mu|} s_{\nu/\mu}(y)\,s_{\mu'}(y)=\delta_{\nu,0}\,s_0(y)\,,
\label{eq-skew-identity}
\end{equation}
for all $y$, where the sum is to be taken over all partitions $\mu$. 

The ring $\Lambda$ of symmetric functions has the structure
of a Hopf algebra, $Symm=(\Lambda,m,\Delta,\iota,\epsilon,\antip)$,
by virtue of the following identification of operators and their action~\cite{FJKW,FJK}:
\begin{itemize}
\item Product $m$:~~~~~~$m(s_\lambda\otimes s_\mu)(\x)=s_\lambda(\x)\cdot s_\mu(\x)=s_\lambda(\x)\,s_\mu(\x)=\sum_{\nu}\,c_{\lambda\mu}^\nu\,s_\nu(\x)$.
\item Unit $\iota$:~~~~~~~~~~~\,$\iota(1)=s_0(\x)$.
\item Coproduct $\Delta$:~~~$s_\nu(\x,\y)=\sum_{\lambda,\mu}\, c_{\lambda\mu}^\nu\,s_\lambda(\x)\,s_\mu(\y)=
\sum_\lambda\ s_\lambda(\x)\otimes s_{\nu/\lambda}(\y)
=\sum_\mu\ s_{\nu/\mu}(\x)\otimes s_\mu(\y)$.
\item Counit $\epsilon$:~~~~~~~~\,$\epsilon(s_\lambda(\x))=\delta_{\lambda0}$.
\item Antipode $\antip$:~~~~~$\antip(s_\lambda(\x))=(-1)^{|\lambda|}\,s_{\lambda'}(\x)$.
\end{itemize}
Here $\x=(x_1,x_2,\ldots)$, $\y=(y_1,y_2,\ldots)$ and $\x,\y=(x_1,x_2,\ldots,y_1,y_2,\ldots)$,
which is sometimes written as $\x+\y$, the addition of two alphabets of indeterminates.
In what follows we tend to favour the use of $\cdot$ rather than $m$ to signify a product.

This particular Hopf algebra $Symm$ has the following properties:
\begin{description}
\item[]Commutativity: $s_\lambda\cdot s_\mu=s_\mu\cdot s_\lambda$ since
$s_\lambda(x)\,s_\mu(x)=s_\mu(x)\,s_\lambda(x)$.
\item[]Cocommutativity: $\Delta(s_\nu)={\sum_\zeta}\ s_\zeta\otimes s_{\nu/\zeta}
={\sum_\zeta}\ s_{\nu/\zeta}\otimes s_\zeta$ since
$s_\nu(x,y)=s_\nu(y,x)$.
\item[]Associativity: $s_\rho\cdot(s_\sigma\cdot s_\tau)=(s_\rho\cdot s_\sigma)\cdot s_\tau$ since
$s_\rho(x)\,\bigl(s_\sigma(x)\,s_\tau(x)\bigr)=\bigl(s_\rho(x)\,s_\sigma(x)\bigr)\,s_\tau(x)$.
\item[]Coassociativity: $(I\otimes\Delta)(\Delta(s_\lambda))=(\Delta\otimes I)(\Delta(s_\lambda))$
since $s_\lambda(x,(y,z))=s_\lambda((x,y),z)$.
\end{description}
In addition we may verify the following general requirements of 
any Hopf algebra:
\begin{description}
\item[]Antipode identity: $\cdot\,(\antip\otimes I)\,\Delta = \iota\, \epsilon
= \cdot\, (I\otimes \antip)\,\Delta$ since, thanks to (\ref{eq-skew-identity}), we have
\begin{eqnarray*}
\cdot\,(\antip\otimes I)\,\Delta(s_\lambda)
&=&\sum_\mu\, \cdot\, (\antip\otimes I)\,(s_\mu\otimes s_{\lambda/\mu})
\ =\ \sum_\mu\, \antip(s_\mu)\cdot s_{\lambda/\mu}\cr
&=&\sum_\mu\,(-1)^{|\mu|}s_{\mu'}\cdot s_{\lambda/\mu}
\ =\ \delta_{\lambda0}\, s_0\ =\ \iota\, \delta_{\lambda0}\ =\
\iota\, \epsilon (s_\lambda)\,.
\end{eqnarray*}
\item[]Counitarity: $\cdot\,(\epsilon\otimes I)\,\Delta=I
= \cdot\,(I\otimes \epsilon)\,\Delta$ since
$$
\cdot\,(\epsilon\otimes I)\,\Delta(s_\lambda)
=\sum_\mu\, \cdot\,(\epsilon\otimes I)\,(s_\mu\otimes s_{\lambda/\mu})
=\sum_\mu\, \cdot\, (\delta_{\mu0}\otimes s_{\lambda/\mu}) = 1 \cdot s_\lambda
=s_\lambda=I(s_\lambda).
$$
\item[]Product and coproduct compatibility: $\Delta\, (\cdot)
=(\cdot\otimes\cdot)\,(\Delta\otimes\Delta)$ since 
\begin{eqnarray*}
\Delta\,(\cdot) (s_\lambda(z)\ s_\mu(w))&=&
\Delta (s_\lambda(z)\ s_\mu(z))\ =\ s_\lambda(x,y)\ s_\mu(x,y)\cr\cr
&=&(\cdot\otimes\cdot)\,(s_\lambda(x,y)\ s_\mu(u,v))
\ =\ (\cdot\otimes\cdot)\,(\Delta\otimes\Delta) (s_\lambda(z)\ s_\mu(w))\,.
\end{eqnarray*}
\end{description}
This last property is the homomorphism property of the coproduct:
\begin{equation}
\Delta(s_\lambda\cdot s_\mu)=\Delta(s_\lambda)\cdot \Delta(s_\mu)
\quad\hbox{or more generally}\quad
\Delta(X\cdot Y)=\Delta(X)\cdot\Delta(Y)\,,
\end{equation}
for any~~$X,Y\in \Lambda$.

At this stage it is convenient to introduce a bilinear scalar
product, $(\cdot\,|\,\cdot)$, on $Symm$. This is defined by
$(s_\lambda\,|\,s_\mu)=\delta_{\lambda\mu}$~\cite{Mac}. With this definition
we have: 
\begin{eqnarray*}
&&(s_\nu\,|\,s_\lambda\cdot s_\mu)=\ds{\sum_\zeta}\ c_{\lambda\mu}^\zeta\ (s_\nu\,|\,s_\zeta)=c_{\lambda\mu}^\nu\,;
\qquad (s_{\nu/\lambda}\,|\,s_\mu)=\ds{\sum_\eta}\ c_{\lambda\eta}^\nu\ (s_\eta\,|\,s_\mu)=c_{\lambda\mu}^\nu\,;\\
&&(\Delta(s_\nu)\,|\,s_\lambda\otimes s_\mu)=
\ds{\sum_{\sigma,\tau}}\ c_{\sigma\tau}^\nu\,(s_\sigma\otimes s_\tau\,|\,s_\lambda\otimes s_\mu)
=\ds{\sum_{\sigma,\tau}}\ c_{\sigma\tau}^\nu\ \delta_{\sigma\lambda}\ \delta_{\tau\mu}
=c_{\lambda\mu}^\nu\,,
\end{eqnarray*}
so that
\begin{equation}
(s_\nu\,|\,s_\lambda\cdot s_\mu)=(s_{\nu/\lambda}\,|\,s_\mu)
\quad\hbox{and}\quad
(s_\nu\,|\,s_\lambda\cdot s_\mu)=(\Delta(s_\nu)\,|\,s_\lambda\otimes s_\mu)\,.
\label{eq-duality}
\end{equation}

For any $X=\sum_\sigma\ a_\sigma\,s_\sigma$ we let
$s_{\lambda\cdot X}=s_\lambda\cdot X=\sum_\sigma\ a_\sigma\ (s_\lambda \cdot s_\sigma)$
and $s_{\lambda/X} =\sum_\sigma\ a_\sigma\ s_{\lambda/\sigma}$.
With this notation, we can extend the first part of (\ref{eq-duality}) so
that for any $X=\sum_\sigma\ a_\sigma\,s_\sigma$ we have 
\begin{equation}
(s_\lambda\,|\,X\cdot s_\mu)=\sum_\sigma\ a_\sigma\ (s_\lambda\,|\,s_\sigma\cdot s_\mu)
=\sum_\sigma\ a_\sigma\ (s_{\lambda/\sigma}\,|\,s_\mu)
=(s_{\lambda/X}\,|\,s_\mu).
\label{eq-lambda-skew-X}
\end{equation}
In addition, it should be noted that for $X=\sum_\sigma\ a_\sigma\,s_\sigma$
we have $a_\sigma=(X\,|\,s_\sigma)$ so that
\begin{equation}
  X=\sum_{\sigma} (X\,|\,s_\sigma)\ s_\sigma
\label{eq-Xsigma}
\end{equation}
Finally, for any $X=\sum_\sigma\ a_\sigma\,s_\sigma$
and $Y=\sum_\sigma\ b_\sigma\,s_\sigma$, we have $X=Y$ if and only 
if $(X\,|\,s_\sigma)=a_\sigma=b_\sigma=(Y\,|\,s_\sigma)$
for all $\sigma$. 
It follows that 
\begin{equation}
s_{(\lambda/\mu)/\nu}=s_{\lambda/(\mu\cdot\nu)}
\label{eq-skew-of-skew}
\end{equation}
since
$(s_{(\lambda/\mu)/\nu}\,|\,s_\sigma)=(s_{\lambda/\mu}\,|\,s_\nu\,s_\sigma)
=(s_{\lambda}\,|\,s_\mu\,s_\nu\,s_\sigma)
=(s_{\lambda}\,|\,s_{\mu\cdot\nu}\,s_\sigma)
=(s_{\lambda/(\mu\cdot\nu)}\,|\,s_\sigma)
$
for all $\sigma$, and
\begin{equation}
s_{(\mu\cdot\nu)/\rho}=\sum_{\sigma,\tau} c^\rho_{\sigma\tau}\ s_{\mu/\sigma}\cdot s_{\nu/\tau},
\label{eq-skew-of-product}
\end{equation} 
since
$\ds
(s_{(\mu\cdot\nu)/\rho}\,|\, s_\lambda)=(s_\mu\cdot s_\nu\,|\,s_\rho\cdot s_\lambda)
=(s_\mu\otimes s_\nu\,|\,\Delta(s_\rho\cdot s_\lambda))
=\sum_{\sigma,\tau} c^\rho_{\sigma\tau}(s_\mu\otimes s_\nu\,|\,s_\sigma\otimes s_\tau\ \Delta(s_\lambda))$\\
$\ds =\sum_{\sigma,\tau} c^\rho_{\sigma\tau}(s_{\mu/\sigma}\otimes s_{\nu/\tau}\,|\,\Delta(s_\lambda))
=\sum_{\sigma,\tau} c^\rho_{\sigma\tau}(s_{\mu/\sigma}\cdot s_{\nu/\tau}\,|\,s_\lambda)
$ for all $\lambda$.

\section{Characters of the classical groups}

Let $M(m,n)$ be the set of all $m\times n$ matrices over $\C$. Then
the classical groups under consideration here are: 
\begin{eqnarray*}
GL(n)&=&\{X\in M(n,n) \,|\, \det X \neq 0\}\,;\cr
O(n)&=&\{X\in GL(n) \,|\, X\,G_n\,X^t=G_n\}\quad\hbox{with $G_n^t=G_n$}\,;\cr
Sp(n)&=&\{X\in GL(n) \,|\, X\,J_n\,X^t=J_n\}\quad\hbox{with $J_n^t=-J_n$}\,.
\end{eqnarray*}
It might be noted that for $n=2k+1$ the matrix $J_n$ is necessarily singular,
and may be chosen~\cite{Pro} so that:
\begin{equation}
Sp(2k+1) =
\left[\begin{array}{cccc}
 &Sp(2k)& &M(2k,1)\cr
 0&\cdots&0&GL(1)\cr
\end{array}\right]\,.
\end{equation}
Thus $Sp(2k+1)$ is not semisimple. Nor is it reductive.

The eigenvalues of an arbitrary group element $X$ may be parametrised
as follows~\cite{Li50}:
\begin{description}
\item[\qquad] $GL(n):~~x_1,x_2,\ldots,x_n$~~with~~$x_1x_2\cdots x_n\neq0$.
\item[\qquad] $SL(n):~~x_1,x_2,\ldots,x_n$~~with~~$x_1x_2\cdots x_n=1$.
\item[\qquad] $SO(2k+1):~~x_1,x_2,\ldots,x_k,\ov{x}_1,\ov{x}_2,\ldots,\ov{x}_k,1$.
\item[\qquad] $O(2k+1)\backslash SO(2k+1):
~~x_1,x_2,\ldots,x_k,\ov{x}_1,\ov{x}_2,\ldots,\ov{x}_k,-1$.
\item[\qquad] $Sp(2k):~~x_1,x_2,\ldots,x_k,\ov{x}_1,\ov{x}_2,\ldots,\ov{x}_k$.
\item[\qquad] $SO(2k):~~x_1,x_2,\ldots,x_k,\ov{x}_1,\ov{x}_2,\ldots,\ov{x}_k$.
\item[\qquad] $O(2k)\backslash SO(2k):
~~x_1,x_2,\ldots,x_{k-1},\ov{x}_1,\ov{x}_2,\ldots,\ov{x}_{k-1},1,-1$.
\item[\qquad] $Sp(2k+1):~~x_1,x_2,\ldots,x_k,\ov{x}_1,\ov{x}_2,\ldots,\ov{x}_k,x_{2k+1}$.
\end{description}
where $\ov{x}_i=x_i^{-1}$ for all $i$.

Each covariant tensor irreducible representation, $V_{GL(n)}^\lambda$, 
of $GL(n)$ is specified by a partition $\lambda$ of length $\ell(\lambda)\leq n$. 
Let $X\in GL(n)$ have eigenvalues $(x_1,x_2,\ldots,x_n)$
and let $\rho=(n-1,n-2,\ldots,1,0)$. 
Then the character of this irreducible representation is given by~\cite{Li50,Mac}:
\begin{equation}
  \ch V_{GL(n)}^\lambda=\frac{a_{\lambda+\rho}(x)}{a_\rho(x)}
  =\frac
   {\left|\,x_i^{\lambda_j+n-j}\,\right|}{\left|\,x_i^{n-j}\,\right|}=s_\lambda(x)\,.
\label{eq-charGL}
\end{equation}
Thanks to the stability property of Schur functions with respect to the
number $n$ of indeterminates, we may define the corresponding universal character 
of $GL(n)$ by
\begin{equation}
  \ch V_{GL}^\lambda=\{\lambda\}(x)=s_\lambda(x)\,,
\end{equation} 
where $\x=(x_1,x_2,\ldots)$. For each finite $n$ the characters $\ch V_{GL(n)}^\lambda$
are recovered from the universal characters $\ch V_{GL}^\lambda$ merely by setting
$\x=(x_1,x_2,\ldots,x_n,0,0,\ldots,0)$.

In a similar way, there exist covariant tensor irreducible representation, $V_{O(n)}^\lambda$
and $V_{Sp(n)}^\lambda$, of $O(n)$ and $Sp(n)$, respectively. The corresponding 
characters $\ch V_{O(n)}^\lambda$ and $\ch V_{Sp(n)}^\lambda$ may each be defined
in terms of determinants. More important, from our point of view, is that
there exist corresponding universal characters~\cite{Ki89,Ko} denoted by
\begin{equation}
  \ch V_{O}^\lambda=[\lambda](x)\quad\hbox{and}\quad \ch V_{Sp}^\lambda=\la\lambda\ra(x)\,,
\label{eq-universal-OSp}
\end{equation} 
with $\x=(x_1,x_2,\ldots)$ arbitrary. These are universal in the sense that for any finite $n$
the characters $\ch V_{O(n)}^\lambda$ and $\ch V_{Sp(n)}^\lambda$ are obtained by
specialising $\x$ to $(x_1,x_2,\ldots,x_n,0,0,\ldots,0)$
with $x_1,x_2,\ldots,x_n$ restricted to the eigenvalues of the
appropriate group elements parametrised as above.

The universal characters (\ref{eq-universal-OSp}) are themselves defined by means 
of the generating functions~\cite{Li50}:
\begin{eqnarray}
\prod_{i,a} (1-x_iy_a)^{-1}\ \prod_{a\leq b} (1-y_ay_b)
&=& \sum_\lambda [\lambda](\x)\ \{\lambda\}(\y)\,;\label{eq-JC}\\
\prod_{i,a} (1-x_iy_a)^{-1}\ \prod_{a< b} (1-y_ay_b)
&=& \sum_\lambda \la\lambda\ra(\x)\ \{\lambda\}(\y)\,.\label{eq-JA}
\end{eqnarray}
This leads to~\cite{Li50,Ki75,Ki89,Ko}:
\begin{thm} The universal characters $\ch V_G^\lambda$ of 
the orthogonal and symplectic groups are given by
\begin{eqnarray}
[\lambda](\x)&=&\{\lambda/C\}(\x)=s_{\lambda/C}(\x)\quad\hbox{where}\quad C(\x)=\prod_{i\leq j}(1-x_ix_j)\,;
\label{eq-O-lambda}\\
\la\lambda\ra(\x)&=&\{\lambda/A\}(\x)=s_{\lambda/A}(\x)\quad\hbox{where}\quad A(\x)=\prod_{i<j}(1-x_ix_j)\,,
\label{eq-Sp-lambda}
\end{eqnarray}
respectively.
\end{thm}

\noindent{\bf Proof:}\ \
For $O(n)$, it follows from (\ref{eq-JC}) that 
the character $[\lambda](x)$ is the coefficient of 
$s_\lambda(y)$ in $J(\x,\y)\, C(\y)$. 
Hence, from (\ref{eq-J}), (\ref{eq-lambda-skew-X})
and (\ref{eq-Xsigma}) 
we have
\begin{eqnarray*}
[\lambda](x)&=&\big( \prod_{i,a} (1-x_iy_a)^{-1}\ \prod_{a\leq b} (1-y_ay_b)\ \big|\ \{\lambda\}(y)\ \big)\\
&=&\big( \sum_\sigma \{\sigma\}(x)\ \{\sigma\}(y)\ C(y)\ \big|\ \{\lambda\}(y)\ \big)\\
&=&\sum_\sigma \{\sigma\}(x)\ \big( \{\sigma\}(y)\cdot C(y)\ \big|\ \{\lambda\}(y)\ \big)\\
&=&\sum_\sigma \{\sigma\}(x)\ \big( \{\sigma\}(y)\ \big|\ \{\lambda/C\}(y)\ \big)=\{\lambda/C\}(x).
\end{eqnarray*}
For $Sp(n)$ one simply replaces
$a\leq b$ by $a<b$, and $C(y)$ by $A(y)$. This gives:
$\la\lambda\ra(x)=\{\lambda/A\}(x)$, as required. 
\qed

\section{Branching rules}

In order to allow the possibility of extending results to
a wider class of subgroups of the general linear group, we
consider any subgroup $H(n)$ of $GL(n)$, whose
group elements $X\in H(n)\subset GL(n)$ have eigenvalues
$\x=(x_1,x_2,\ldots,x_n)$. We assume, just as in the case of $O(n)$ and $Sp(n)$,
that there exist irreducible representations $V_{H(n)}^\lambda$ 
of $H(n)$, specified by partitions $\lambda$, with characters
$\ch V_{H(n)}^\lambda$ that may be determined by
specialising from $\x=(x_1,x_2,\ldots)$ to $\x=(x_1,x_2,\ldots,x_n,0,\ldots,0)$,
with appropriate $x_1,x_2,\ldots,x_n$,
the universal characters
\begin{equation}
\ch V_H^\lambda = \lb\lambda\rb(\x)\,.
\label{eq-Hchar}
\end{equation}
If the embedding $H(n)\subset GL(n)$ is such that there exists two 
series of Schur functions $S(x)$ and $T(x)$ with the properties:
\begin{equation}
\lb\lambda\rb(x)=\{\lambda/S\}(x)
\quad\hbox{with}\quad S(x)\ T(x) =1
\label{eq-S}
\end{equation}
then 
\begin{equation}
\{\lambda\}(x)=\lb\lambda/T\rb(x)\,.
\label{eq-T}
\end{equation}
As a result we immediately have the branching rule:
\begin{equation}
GL(n)\supset H(n):\qquad \{\lambda\}\rightarrow\lb\lambda/T\rb\,.
\label{eq-GLbrH}
\end{equation}

Applying this to $O(n)$ and $Sp(n)$, we immediately have~\cite{Li50,Ki75}:
\begin{thm} The branching rules for the decomposition of
representations of $GL(n)$ under restriction to the subgroups $O(n)$
and $Sp(n)$ take the form:
\begin{eqnarray}
GL(n)\supset O(n):\ && \{\lambda\}\rightarrow[\lambda/D]
\quad\hbox{with}\quad D=C^{-1}=\prod_{i\leq j} (1-x_ix_j)^{-1}\,;\\
GL(n)\supset Sp(n):\ && \{\lambda\}\rightarrow\la\lambda/B\ra
\quad\hbox{with}\quad B=A^{-1}=\prod_{i<j} (1-x_ix_j)^{-1}\,.
\end{eqnarray}
\end{thm}

It is well known that~\cite{Ki75,BKW}
\begin{eqnarray*}
B&=&\{0\}+\{1^2\}+\{2^2\}+\{1^4\}+\{3^2\}+\{2^21^2\}+\{1^6\}+\cdots\\
D&=&\{0\}+\{2\}+\{4\}+\{2^2\}+\{6\}+\{42\}+\{2^3\}+\cdots
\end{eqnarray*}
where $D$ involves partitions all of whose parts are even, 
and $B$ the conjugate of such partitions.
The branching rules obtained using this identification
of $D$ and $B$ are exemplified in Table~\ref{tab-branching}.

\begin{center}
\begin{table}[h]
\caption{\label{tab-branching} Branching rule examples.}
\centering
\begin{tabular}{@{}*{7}{l}}
\br
$GL(n)\supset O(n)$:&$\{\lambda\}\rightarrow[\lambda/D]$\\
\mr
&$\{4\}\rightarrow [4]+[2]+[0]$\\
&$\{1^4\}\rightarrow [1^4]$\\
&$\{2^21^2\}\rightarrow [2^21^2]+[21^2]+[1^2]$\\
\mr
$GL(n)\supset Sp(n)$:& $\{\lambda\}\rightarrow\la\lambda/B\ra$\\
\mr
&$\{4\}\rightarrow \la4\ra$\\
&$\{1^4\}\rightarrow \la1^4\ra+\la1^2\ra+\la0\ra$\\
&$\{2^21^2\}\rightarrow \la2^21^2\ra+\la2^2\ra+\la21^2\ra+\la1^4\ra+2\,\la1^2\ra+\la0\ra$\\
\br
\end{tabular}
\end{table}
\end{center}

\section{Tensor products}

We return to the general case of the subgroup $H(n)$ of $GL(n)$ with
universal characters defined by $\lb\lambda\rb=\{\lambda/S\}$.
Let the coproduct of $T=S^{-1}$ take the form
\begin{equation}
\Delta(T)=(T\otimes T)\cdot \Delta''(T)
\qquad\hbox{with}\qquad
\Delta''(T)=\sum_{\sigma,\tau}
b_{\sigma\tau}^T\ \{\sigma\}\otimes \{\tau\}\,,
\label{eq-coprod-T}
\end{equation}
for some, as yet undetermined, coefficients $b_{\sigma\tau}^T$.
Then, we have
\begin{thm}
The decomposition of products of universal characters
of $H(n)$ takes the form:
\begin{equation}
\lb\lambda\rb\cdot\lb\mu\rb =\sum_{\sigma,\tau} b_{\sigma\tau}^T\
\lb(\lambda/\sigma)\cdot(\mu/\tau\rb.
\label{eq-tensor-prod}
\end{equation}
\end{thm}

\noindent{\bf Proof:}\ \
Note that from (\ref{eq-S}) and (\ref{eq-T}) 
\begin{equation}
\quad\lb\lambda\rb\cdot\lb\mu\rb =\{\lambda/S\}\cdot \{\mu/S\}
=\lb((\lambda/S)\cdot(\mu/S))/T\rb\,.
\end{equation}
This implies that the multiplicity of $\lb\nu\rb$ in
$\lb\lambda\rb\cdot\lb\mu\rb$ is the same as the multiplicity of
$\{\nu\}$ in $\{((\lambda/S)\cdot(\mu/S))/T\}$, that is:
\begin{eqnarray*}
(\{((\lambda/S)\cdot(\mu/S))/T\}\ |\ \{\nu\})
&=&(\{(\lambda/S)\cdot(\mu/S)\}\ |\ T\cdot\{\nu\})\\
&=&(\{\lambda/S\}\otimes\{\mu/S\}\ |\ \Delta(T\cdot\{\nu\}))\\
&=&(\{\lambda/S\}\otimes\{\mu/S\}\ |\ (T\otimes T)\cdot\Delta''(T)\cdot\Delta(\{\nu\}))\\
&=&(\{\lambda/ST\}\otimes\{\mu/ST\}\ |\
\Delta''(T)\cdot\Delta(\{\nu\}))\\
&=&(\{\lambda\}\otimes\{\mu\}\ |\ \sum_{\sigma,\tau}
b_{\sigma\tau}^T\
\{\sigma\}\otimes\{\tau\}\cdot\Delta(\{\nu\}))\\
&=&\sum_{\sigma,\tau} b_{\sigma\tau}^T\
(\{\lambda/\sigma\}\otimes\{\mu/\tau\}\ |\ \Delta(\{\nu\}))\\
&=&\sum_{\sigma,\tau} b_{\sigma\tau}^T
(\{\lambda/\sigma\}\cdot\{\mu/\sigma\}\ |\ \{\nu\})\,,
\end{eqnarray*} 
from which the required formula (\ref{eq-tensor-prod})
follows
\qed

In order to apply this to any given subgroup $H(n)$ of $GL(n)$
it is necessary to evaluate the coefficients $b_{\sigma\tau}^T$
appearing in the coproduct of $T$. In the case of $O(n)$
and $Sp(n)$ we have $T=D$ and $T=B$, respectively, for which
we have the coproduct expansions
\begin{eqnarray}
\Delta(D)&=&(D\otimes D)\cdot\Delta''(D)\quad\hbox{with}\quad \Delta''(D)
=\sum_\sigma \{\sigma\}\otimes\{\sigma\}\,;
\label{eq-coproduct-D}\\
\Delta(B)&=&(B\otimes B)\cdot\Delta''(B)\quad\hbox{with}\quad \Delta''(B)
=\sum_\sigma \{\sigma\}\otimes\{\sigma\}\,.
\label{eq-coproduct-B}
\end{eqnarray}
This can be seen by noting that:
\begin{eqnarray*}
D(\x,\y)&=&\prod_{i\leq j} (1- x_ix_j)^{-1}\prod_{i,a} (1-x_iy_a)^{-1}\prod_{a\leq b} (1-y_ay_b)^{-1}\\
&=&D(\x)\ \sum_\sigma \{\sigma\}(\x)\,\{\sigma\}(\y)\ D(\y)\,;\\
B(\x,\y)&=&\prod_{i<j} (1- x_ix_j)^{-1}\prod_{i,a} (1-x_iy_a)^{-1}\prod_{a<b} (1-y_ay_b)^{-1}\\
&=&B(\x)\ \sum_\sigma \{\sigma\}(\x)\,\{\sigma\}(\y)\ B(\y)\,.
\end{eqnarray*}

Now we are in a position to apply (\ref{eq-tensor-prod}) to
the case of $O(n)$ and $Sp(n)$. We find~\cite{New,Li58,BKW,Ki89,Ko}:

\begin{thm}
\label{thm-tensor-prod-OSp}
The tensor product decomposition rules for universal characters 
of $O(n)$ and $Sp(n)$ take the form:
\begin{equation}
[\lambda]\,\cdot\,[\mu]=\sum_\sigma\
[(\lambda/\sigma)\cdot(\mu/\sigma)]
\qquad\hbox{and}\qquad
\la\lambda\ra\,\cdot\,\la\mu\ra=\sum_\sigma\
\la(\lambda/\sigma)\cdot(\mu/\sigma)\ra\,.
\label{eq-prod-OSp}
\end{equation}
\end{thm}

\noindent{\bf Proof:}\ \ 
In the case of $O(n)$ it is merely necessary to note that $T=D$
and $b_{\sigma\tau}^D=\delta_{\sigma,\tau}$ from (\ref{eq-coproduct-D}).
Using this in (\ref{eq-tensor-prod}) gives
\begin{equation}
[\lambda]\cdot[\mu]=\sum_{\sigma,\tau}
b_{\sigma\tau}^D\ [(\lambda/\sigma)\cdot(\mu/\tau)]
=\sum_\sigma\ [(\lambda/\sigma)\cdot(\mu/\sigma)]\,.
\end{equation}
Similarly, in the case $Sp(n)$
we have $T=B$ and $b_{\sigma\tau}^B=\delta_{\sigma,\tau}$ from 
(\ref{eq-coproduct-B}). Again using this in (\ref{eq-tensor-prod})
gives
\begin{equation}
\la\lambda\ra\cdot\la\mu\ra=\sum_{\sigma,\tau}
b_{\sigma\tau}^B\ \la(\lambda/\sigma)\cdot(\mu/\tau)\ra
=\sum_\sigma\ \la(\lambda/\sigma)\cdot(\mu/\sigma)\ra\,.
\end{equation}
\qed

These rules are exemplified in Table~\ref{tab-products}, along
with an example for $GL(n)$ that is included for comparative purposes.
Remarkably, the universal tensor product rules for $O(n)$ and $Sp(n)$ are
identical. However, it is important to note that for finite $n$, 
when $\x=(x_1,x_2,\ldots)$ is suitably specialised, modification rules
given elsewhere~\cite{Ki71,BKW} distinguish them.

\begin{center}
\begin{table}[h]
\caption{\label{tab-products} Tensor product decompositions.}
\centering
\begin{tabular}{@{}*{7}{l}}
\br
$GL(n)$:&$\ \{2^2\}\cdot\{21\}$
&$=\{43\}+\{421\}+\{3^21\}+\{32^2\}+\{321^2\}+\{2^31\}$\\ 
\mr
$O(n)$:&$\ [2^2]\cdot[21]$
&$=[43]+[421]+[3^21]+[32^2]+[321^2]+[2^31]$\\
&&$\ \ +[41]+2[32]+2[31^2]+2[2^21]+[21^3]$\\
&&$\ \ +[3]+2[21]+[1^3]+[1]$.\\ 
\mr
$Sp(n)$:&$\ \la2^2\ra\cdot\la21\ra$
&$=\la43\ra+\la421\ra+\la3^21\ra+\la32^2\ra+\la321^2\ra+\la2^31\ra$\\
&&$\ \ +\la41\ra+2\la32\ra+2\la31^2\ra+2\la2^21\ra+\la21^3\ra$\\
&&$\ \ +\la3\ra+2\la21\ra+\la1^3\ra+\la1\ra$.\\
\br
\end{tabular}
\end{table}
\end{center}

\section{Classical group character rings}

Following an approach described in~\cite{FJK}, 
the universal characters $\{\lambda\}$, $[\lambda]$ and 
$\la\lambda\ra$, 
of the general linear, othogonal and symplectic
groups are linked to each other 
by means of the following identities:
\begin{equation}
\{\lambda\}=[\lambda/D]=\la\lambda/B\ra\,;
\quad
\{\lambda/C\}=[\lambda]=\la\lambda/BC\ra\,;
\quad
\{\lambda/A\}=[\lambda/AD]=\la\lambda\ra\,.
\label{eq-bases-gosp} 
\end{equation}
By virtue of these identities each of these sets of characters
$\{\lambda\}$, $[\lambda]$ and $\la\lambda\ra$ forms a basis
of $\Lambda$. However, as we have seen, their product rules within the character
rings $CharGL$, $CharO$ and $CharSp$ are different.
They are tabulated in the first line of Table~\ref{tab-char}.
Moreover, within these same character rings, the coproduct
of the Hopf algebra $Symm$ induces the coproducts
given in the second line of Table~\ref{tab-char}.
These coproduct formulae are just the universal form of the 
branching rules~\cite{Ki75} for the group-subgroup
restrictions $GL(n+m)\supset GL(n)\times GL(m)$, $O(n+m)\supset O(n)\times O(m)$,
and $Sp(n+m)\supset Sp(n)\times Sp(m)$.
To complete the specification of the Hopf algebra structure of $CharGL$, $CharO$ and 
$CharSp$ it is only necessary to specify the unit $\iota$, counit $\epsilon$ and 
antipode $\antip$. These are also given in Table~\ref{tab-char}, where $\cal A$
and $\cal C$ signify the sets of partitions $\alpha$ and $\gamma$ appearing in the
expansions of $A$ and $C$~\cite{Ki75,BKW}.

\begin{center}
\begin{table}[h]
\caption{\label{tab-char} Hopf algebra structure of group character rings~\cite{FJK}}
\centering
\begin{tabular}{@{}*{7}{l}}
\br
$CharGL$&$CharO$&$CharSp$\\
\br
$m(\{\mu\}\otimes\{\nu\})=\{\mu\}\cdot\{\nu\}$
&
$m([\mu]\otimes[\nu])=\sum_\zeta [(\mu/\zeta)\cdot(\nu/\zeta)]$
&
$m(\la\mu\ra\otimes\la\nu\ra)=\sum_\zeta \la(\mu/\zeta)\cdot(\nu/\zeta)\ra$\\ \\
$\Delta(\{\lambda\})=\sum_{\zeta} \{\lambda/\zeta\}\otimes\{\zeta\}$
&
$\Delta([\lambda])=\sum_{\zeta} [\lambda/\zeta]\otimes[\zeta/D]$
&
$\Delta(\la\lambda\ra)=\sum_{\zeta} \la\lambda/\zeta\ra\otimes\la\zeta/B\ra$\\ \\
$\iota(1)=\{0\}$&$\iota(1)=[0]$&$\iota(1)=\la0\ra$\\ \\
$\epsilon(\{\lambda\})=\delta_{\lambda,0}$
&
$\epsilon([\lambda])=\sum_{\gamma\in{\cal C}} (-1)^{|\gamma|/2}\,\delta_{\lambda,\gamma}$
&
$\epsilon(\la\lambda\ra)=\sum_{\alpha\in{\cal A}} (-1)^{|\alpha|/2}\,\delta_{\lambda,\alpha}$\\ \\
$\antip(\{\lambda\})=(-1)^{|\lambda|}\{\lambda'\}$
&
$\antip([\lambda])=(-1)^{|\lambda|}[\lambda'/AD]$
&
$\antip(\la\lambda\ra)=(-1)^{|\lambda|}\la\lambda'/CB\ra$\\
\br
\end{tabular}
\end{table}
\end{center}

\section{Conclusions}

Universal characters of irreducible 
representations of the classical groups have been identified.
These have been expressed in terms of Schur functions, and
Hopf algebra manipulations have allowed us to calculate
branching rules and tensor product decompositions.
The analysis covers covariant tensor 
representations of $GL(n)$, $O(n)$ and $Sp(n)$,
and may be extended to other subgroups of $GL(n)$~\cite{FJKW}.
Such an extension using higher rank 
invariants, leads in some cases to finite subgroups. 

It should be stressed that, for any finite $n$, modification rules 
are needed to interpret the results~\cite{New,Ki71,BKW}. 
In addition, for $n$ odd, the subgroup $Sp(n)$ of $GL(n)$ is neither
semisimple nor reductive. However, for $n=2k+1$ the results remain 
valid if each character $\la\lambda\ra$ is interpreted as the character of a 
representation $V_{Sp(2k+1)}^\lambda$ that is no longer irreducible
but is indecomposable~\cite{Pro}.

\end{document}